\documentclass[12pt]{amsart}
\usepackage{amsfonts}
\usepackage{graphicx}
\usepackage{amsmath}
\usepackage{amssymb}
\input{amssym.def}
\newtheorem{theorem}{Theorem}

\theoremstyle{remark}

\newcommand{\re}{\text{\rm Re }}
\newcommand{\im}{\text{\rm Im }}

\input epsf
\begin{document}
\title[Harmonic measures]{Asymptotic ratio of harmonic measures of slit sides}
\author[D.~Prokhorov, D.~Ukrainskii]{Dmitri Prokhorov and Dmitrii Ukrainskii}

\address{D.~Prokhorov, D.~Ukrainskii: Department of Mathematics and Mechanics, Saratov State
University, Saratov 410012, Russia} \email{ProkhorovDV@info.sgu.ru}
\email{D.v.ukrainskiy@gmail.com}
\thanks{D.Prokhorov have been supported by the Russian Ministry of Education and Science (project
1.1520.2014k). D.Ukrainskii have been supported by Russian/Turkish grant RFBR/T\"UBITAK
\#14-01-91370} 

\begin{abstract}
The article is devoted to the geometry of solutions to the chordal L\"owner equation which is based
on the comparison of singular solutions and harmonic measures for the sides of a slit in the upper
half-plane generated by a driving term. An asymptotic ratio for harmonic measures of slit sides is
found for a slit which is tangential to a straight line under a given angle, and for a slit with
high order tangency to a circular arc tangential to the real axis.
\end{abstract}
\maketitle

\section{Introduction}

The famous L\"owner differential equation have been introduced in 1923 \cite{Lowner} and was aimed
to give a parametric representation of slit domains. In this article we describe an asymptotic
behavior of singular solutions and harmonic measures for the sides of a slit in domains generated
by a driving term of the L\"owner equation.

The chordal version of the L\"owner equation deals with the upper half-plane $\mathbb H=\{z: \im
z>0\}$, $\mathbb R=\partial \mathbb H$, and functions $f(z,t)$ normalized near infinity by $$
f(z,t)=z+\frac{2t}{z}+O\left(\frac{1}{z^2}\right)$$ which solve the chordal L\"owner differential
equation
\begin{equation}
\frac{df(z,t)}{dt}=\frac{2}{f(z,t)-\lambda(t)}, \quad f(z,0)\equiv z, \quad t\geq0, \label{Low}
\end{equation}
and map subdomains of $\mathbb H$ onto $\mathbb H$. Here $\lambda(t)$ is a real-valued
continuous driving term.

Let $\gamma_t:=\gamma[0,t]=\{\gamma(x): 0\leq x\leq t\}$ be a simple continuous curve in $\mathbb
H\cup\{0\}$ with endpoints $\gamma(0)=0$ and $\gamma(t)$, $0\leq t\leq T$. Then there is a unique
map $f(z,t): \mathbb H\setminus\gamma_t\to\mathbb H$ satisfying the chordal L\"owner equation
(\ref{Low}) with $\lambda(t)$ uniquely determined by $\gamma[0,t]$. The function $f(z,t)$ can be
extended continuously to $\mathbb R\cup\gamma(t)$, and $f(\gamma(t),t)=\lambda(t)$. The value $t$
is called the {\it half-plane capacity} of the curve $\gamma_t$, $t=\text{hcap}(\gamma_t)$, see,
e.g., \cite{Lind}.

We say that $\gamma_t\in C^n$, $n\in\mathbb N$, on $[0,S]$ if, for the arc-length parameter $s$ of
$\gamma_t$, $\gamma(t(s))$ has a continuous derivative $\gamma^{(n)}$ in $s$ on $[0,S]$, $t(S)=T$.
All the derivatives $\gamma^{(k)}$, $1\leq k\leq n$, at $s=0$ are understood as one-side
derivatives. A curve $\gamma_t\in C^n$, $\gamma(0)=0$, is said to have at least $n$-order tangency
with a ray $I_{\theta}=\{e^{i\theta}s: s\geq0\}$, $\theta\in\mathbb R$, at $s=0$ if
$$\gamma(t(s))=e^{i\theta}s+o(s^n),\;\;\;s\to+0.$$ Two curves $\gamma[0,s]\in C^n$ and
$\Gamma[0,s]\in C^n$ are said to have at least $n$-order tangency at $s=0$ if derivatives
$\gamma^{(k)}$ and $\Gamma^{(k)}$ in $s$ at $s=0$ coincide, $0\leq k\leq n$.

The extended function $f(z,t)$ maps $\gamma_t$ onto a segment $I=I(t)=[f_2(0,t),f_1(0,t)]$ while
$f(\mathbb R)=\mathbb R\setminus I$. The function $f_1(0,t)$ is the maximal singular solution to
the chordal L\"owner equation (\ref{Low}), and $f_2(0,t)$ is the minimal singular solution to
(\ref{Low}). Both of these solutions correspond to the singular point $f(0,0)=0$ of equation
(\ref{Low}).

The curve $\gamma_t$ has two sides $\gamma_{1t}$ and $\gamma_{2t}$ which define different prime
ends at the same points, except for its tip. We say that $\gamma_{1t}$ is the {\it left} side of
$\gamma_t$ if going along the boundary of the domain $\mathbb H\setminus\gamma_t$ and moving along
$\mathbb R$ from $(-\infty)$ to 0, we first meet the side $\gamma_{1t}$ and then $\gamma_{2t}$. In
this case, $\gamma_{2t}$ is called the {\it right} side of $\gamma_t$. The two parts
$[f_2(0,t),\lambda(t)]$ and $[\lambda(t),f_1(0,t)]$ of segment $I(t)$ are the images of the two
sides $\gamma_{1t}$ and $\gamma_{2t}$ of $\gamma_t$ under $f(z,t)$, respectively.

The harmonic measures $\omega(f^{-1}(i,t);\gamma_{kt},\mathbb H\setminus\gamma_t)$ of $\gamma_{kt}$
at $f^{-1}(i,t)$ with respect to $\mathbb H\setminus\gamma_t$ are defined by the functions
$\omega_k$ which are harmonic on $\mathbb H\setminus\gamma_t$ and continuously extended on its
closure except for the endpoints of $\gamma_t$, $\omega_k|_{\gamma_{kt}}=1$, $\omega_k|_{\mathbb
R\cup(\gamma_t\setminus\gamma_{kt})}=0$, $k=1,2$, see, e.g., [3, \S3.6]. Denote
$$m_k(t):=\omega(f^{-1}(i,t);\gamma_{kt},\mathbb H\setminus\gamma_t),\;\;\;k=1,2.$$

In Section 2, we prove the following theorem.

\begin{theorem}
Let $\gamma_t\in C^4$, $\gamma(0)=0$, $\im\gamma(t)>0$ for $t>0$, have at least 4-order tangency at
the origin to the straight line under the angle ${\pi\over2}(1-\beta)$, $-1<\beta<1$, to the real
axis $\mathbb R$, and let $f(z,t)$ map $\mathbb H\setminus\gamma_t$ onto $\mathbb H$ and solve the
chordal L\"owner equation (\ref{Low}). Then
$$\lim_{t\to+0}\frac{m_1(t)}{m_2(t)}=\frac{1+\beta}{1-\beta},$$ where
$$m_k(t):=\omega(f^{-1}(i,t);\gamma_{kt},\mathbb H\setminus\gamma_t),\;\;\;k=1,2,$$ $\gamma_{1t}$
is the left side of $\gamma_t$, and $\gamma_{2t}$ is the right side of $\gamma_t$.
\end{theorem}

The most important argument in the proof is the comparison of asymptotic parametric representations
of $\gamma_t$ in $t$ and $s$ at $s=0$. This approach can be compared with the result by Earle and
Epstein \cite{Earle}.

In Section 3, we solve a similar problem for a curve $\gamma_t$ which has at least 6-order tangency
with a circular arc in $\mathbb H\cup\{0\}$ tangential to $\mathbb R$ at the origin. Since a
scaling time change $t\to\alpha^2t$ in the L\"owner equation (\ref{Low}) is accompanied by changing
$\lambda(t)\to{1\over\alpha}\lambda(\alpha^2t)$ and $f(z,t)\to{1\over\alpha}f(\alpha z,\alpha^2t)$,
$\alpha\in\mathbb R$, we can assume without loss of generality that the circular arc is of radius
1, and the argument of its points is increasing when going from 0. We prove the following theorem.

\begin{theorem}
Let $\gamma_t\in C^6$, $\gamma(0)=0$, $\im\gamma(t)>0$ and $\re\gamma(t)>0$ for $t>0$, have at
least 6-order tangency at the origin to the circular arc of radius 1 centered at $i$, and let
$f(z,t)$ map $\mathbb H\setminus\gamma_t$ onto $\mathbb H$ and solve the chordal L\"owner equation
(\ref{Low}). Then
$$\lim_{t\to+0}\frac{M^2_1(t)}{M_2(t)}=2\pi,$$ where
$$M_k(t):=\omega(f^{-1}(i,t);\gamma_{kt},\mathbb H\setminus\gamma_t),\;\;\;k=1,2,$$ $\gamma_{1t}$
is the left side of $\gamma_t$, and $\gamma_{2t}$ is the right side of $\gamma_t$.
\end{theorem}

\section{Proof of Theorem 1}

{\it Proof of Theorem 1.} For $\beta=0$, Theorem 1 has been proved \cite{Zakharov}. The cases
$\beta>0$ and $\beta<0$ are symmetric to each other, and we will stop only on $\beta>0$.

The L\"owner equation (\ref{Low}) can be integrated in quadratures in particular cases
\cite{Kager}. For example, if $\lambda(t)=c\sqrt t$, $c\geq0$, then a solution $f_c(z,t)$ to
equation (\ref{Low}) maps $\mathbb H\setminus\gamma_t$ onto $\mathbb H$ where $\gamma_t$ is
parameterized as
$$\gamma[0,t]=\{z=B\sqrt x: 0\leq x\leq t\},$$ with $B=B(c)=|B(c)|e^{i\theta(c)}$,
$$|B(c)|=2\left(\frac{\sqrt{c^2+16}+c}{\sqrt{c^2+16}-c}\right)^{\frac{c}{2\sqrt{c^2+16}}},$$
$$\theta(c)=\frac{\pi}{2}\left(1-\frac{c}{\sqrt{c^2+16}}\right).$$

Suppose that a $C^4$-slit $\gamma_t$ satisfies the conditions of Theorem 1. Then there exists a
driving function $\lambda(t)\in\text{Lip}({1\over2})$ such that a solution $w=f(z,t)$ to equation
(\ref{Low}) maps $\mathbb H\setminus\gamma_t$ onto $\mathbb H$. For the arc-length parameter $s$,
$\gamma(t(s))$ is represented as
\begin{equation}
\gamma(t(s))=e^{i\theta}s+o(s^4),\;\;\;s\to+0. \label{arc}
\end{equation}
Denote $$I_{\theta}(t)=\{xe^{i\theta}: 0\leq x\leq t\},\;\;0<\theta<\frac{\pi}{2},\;\;t>0.$$

There is $c>0$ such that
\begin{equation}
\beta=\frac{c}{\sqrt{c^2+16}} \label{beta}
\end{equation}
for which $\theta=\theta(c)={\pi\over2}(1-\beta)$. Then $f_c(z,\tau)$ maps $\mathbb H\setminus
I_{\theta(c)}(|B(c)|\sqrt{\tau})$ onto $\mathbb H$. The length $\sigma(\tau)$ of
$I_{\theta(c)}(|B(c)|\sqrt{\tau})$ and the half-plane capacity $\tau$ of
$I_{\theta(c)}(|B(c)|\sqrt{\tau})$ are related by
\begin{equation}
\sigma(\tau)=|B(c)|\sqrt{\tau},\;\;\;\tau>0. \label{tau}
\end{equation}

Let $s$ denote the length of $\gamma_{t(s)}$, and let $\sigma$ denote the length of projection of
$\gamma_{t(s)}$ onto $I_{\theta(c)}(T)$ for $T$ large enough. There is a $C^4$-dependence
$s=s(\sigma)$, $$s(0)=0,\;\;s'(0)=1,\;\;s^{(k)}(0)=0,\;k=2,3,4.$$ Therefore,
\begin{equation}
s=\sigma+o(\sigma^4),\;\;\sigma\to+0. \label{sigma}
\end{equation}

Asymptotic expansion (\ref{arc}) implies an asymptotic behavior of a distance between $\gamma_t$
and its projection on $I_{\theta}$,
$$\text{dist}(\gamma_{t(s)},I_{\theta(c)}(\sigma(s)))=\sigma+o(\sigma^4(s)),\;\;s\to+0.$$

Lind, Marshall and Rohde \cite{Lind} studied the closeness of half-plane capacities for two curves
which are close together. According to Lemma 4.10 \cite{Lind}, we have that $$
t(s)-\tau(\sigma(s))=o(s^2),\;\;s\to+0,$$ where $\sigma(\tau)$ is given by (\ref{tau}). Hence, due
to (\ref{tau}) and (\ref{sigma}),
$$t(s(\sigma))=\tau(\sigma)+o(s^2(\sigma))=\tau(\sigma)+o(\sigma^2)=|B(c)|^{-2}s^2+o(s^2),\;\;
s\to+0.$$ Take into account (\ref{arc}) and rewrite the last relation in the form
\begin{equation}
\gamma(t)=|B(c)|\sqrt t+\alpha(t)\sqrt t,\;\;\;\lim_{t\to+0}\alpha(t)=0. \label{basic}
\end{equation}

Choose an arbitrary sequence $\{x_n\}$ of positive numbers $x_n$, $x_n\to\infty$ as $n\to\infty$,
and denote
\begin{equation}
z=g_n(w,t):=\sqrt{x_n}f^{-1}\left(\frac{w}{\sqrt{x_n}},\frac{t}{x_n}\right),\;\;\;n=1,2, \dots\;.
\label{seq} \end{equation} The function $z=f^{-1}(w,t)$ maps $\mathbb H$ onto $\mathbb
H\setminus\gamma[0,t]$, $f^{-1}(\lambda(t),t)=\gamma(t)$. So the functions $g_n(w,t)$ map $\mathbb
H$ onto $\mathbb H\setminus\gamma^{(n)}[0,t]$ where
$$\gamma^{(n)}(t)=\sqrt{x_n}\,\gamma\left(\frac{t}{x_n}\right)=
e^{i\theta(c)}|B(c)|\sqrt t+\alpha\left(\frac{t}{x_n}\right)\sqrt t,$$ $$g_n\left(\sqrt{x_n}\;
\lambda\left(\frac{t}{x_n}\right),t\right)=\gamma^{(n)}(t),\;\;\;0<t\leq T.$$

We see that $$\gamma^{(n)}(t)-e^{i\theta(c)}|B(c)|\sqrt t=\alpha\left(\frac{t}{x_n}\right)\sqrt
t\to0,\;\;n\to\infty,$$ and the convergence is uniform with respect to $t\in[0,T]$.

The Rad\'o theorem \cite{Rad}, see also [2, p.60], states that a sequence $\{h_n\}$ of conformal
mappings $h_n$ from the unit disk $\mathbb D$ onto simply connected domains $D_n$ bounded by Jordan
curves $\partial D_n$, $0\in D_n$, $h_n(0)=0$, $h_n'(0)>0$, converges uniformly on the closure of
$\mathbb D$ to $h:\mathbb D\to D$, $\partial D$ is bounded by a Jordan curve, if and only if $D_n$
converges to the kernel $D$ and, for every $\epsilon>0$, there exists $N>0$ such that, for all
$n>N$, there is a one-to-one correspondence $z_n:\partial D_n\to \partial D$,
$|z_n(\zeta)-\zeta|<\epsilon$, $\zeta\in\partial D_n$. Markushevich \cite{Mar} generalized the
Rad\'o theorem to domains with arbitrary boundaries.

Apply the Rad\'o-Markushevich theorem to $g_n\circ p$ with a conformal mapping $p$ from $\mathbb D$
onto $\mathbb H$ and obtain that the sequence $\{g_n(w,t)\}$ converges to $f_c^{-1}(w,t)$ as
$n\to\infty$ uniformly on compact subsets of $\mathbb H\cup\mathbb R$.

Denote by $\Gamma_1[0,\tau(t)]$ the {\it left} side of the segment
$I_{\theta(c)}(|B(c)|\sqrt{\tau})$ and denote by $\Gamma_2[0,\tau(t)]$ the {\it right} side of this
segment. Similarly, denote $\gamma_{1n}[0,t]$ the {\it left} side of $\gamma^{(n)}[0,t]$, and
denote $\gamma_{2n}[0,t]$ the {\it right} side of $\gamma^{(n)}[0,t]$. The functions
$g_n^{-1}(z,t)$ map $\gamma_{1n}[0,t]$ and $\gamma_{2n}[0,t]$ onto segments
$I_{1n}=I_{1n}(t)\subset\mathbb R$ and $I_{2n}=I_{2n}(t)\subset\mathbb R$, respectively. It is
known, see, e.g., \cite{Zakharov}, that slit sides $\Gamma_1[0,\tau(t)]$ and $\Gamma_2[0,\tau(t)]$
are mapped by $f_c(z,t)$ onto $$I_1=I_1(t)=\left[\frac{c-\sqrt{c^2+16}}{2}\sqrt t,c\sqrt
t\right]\;\;\text{and}\;\;I_2=I_2(t)=\left[c\sqrt t,\frac{c+\sqrt{c^2+16}}{2}\sqrt t\right],$$
respectively. The uniform convergence of $g_n$ to $f_c^{-1}$ implies that $I_{1n}(t)$ tend to
$I_1(t)$, and $I_{2n}(t)$ tend to $I_2(t)$ as $n\to\infty$.

Denote by $\gamma'_{1n}[0,t]$ and $\gamma'_{2n}[0,t]$ the {\it left} and the {\it right} sides of
$\gamma[0,{t\over x_n}]$, respectively. The function $f(z,{t\over x_n})$ maps slit sides
$\gamma'_{1n}[0,t]$ and $\gamma'_{2n}[0,t]$ onto segments $I'_{1n}=I'_{1n}(t)\subset\mathbb R$ and
$I'_{2n}=I'_{2n}(t)\subset\mathbb R$, respectively. Compare $I_{kn}(t)$ and $I'_{kn}(t)$ by
(\ref{seq}) and conclude that $I_{kn}(t)=\sqrt{x_n}I'_{kn}(t)$, and so
$$\text{meas}\;I_{kn}(t)=\sqrt{x_n}\,\text{meas}\;I'_{kn}(t),\;\;k=1,2,\;\;n\geq1,\;\;0<t\leq T.$$

The harmonic measure is invariant under conformal transformations. This gives that $$
\frac{m_1({t\over x_n})}{m_2({t\over x_n})}=\frac{\omega(f^{-1}(i,{t\over x_n}),\gamma_1[0,{t\over
x_n}],\mathbb H\setminus\gamma({t\over x_n}))}{\omega(f^{-1}(i,{t\over x_n}),\gamma_2[0,{t\over
x_n}],\mathbb H\setminus\gamma({t\over x_n}))}=\frac{\omega(i,I'_{1n}(t),\mathbb
H)}{\omega(i,I'_{2n}(t),\mathbb H)}.$$

For $k=1,2$, $n\geq1$, the harmonic measure $\omega(i;I'_{kn}(t),\mathbb H)$ of $I'_{kn}(t)$ at $i$
with respect to $\mathbb H$ equals the angle divided over $\pi$ under which the segment
$I'_{kn}(t)$ is seen from the point $i$. Similarly, the harmonic measure
$\omega(i;I_{kn}(t),\mathbb H)$ of $I_{kn}(t)$ at $i$ with respect to $\mathbb H$ equals the angle
divided over $\pi$ under which the segment $I_{kn}(t)$ is seen from the point $i$, see, e.g., [2,
p.334]. This shows that the last term in the chain of equalities has a limit as $n\to\infty$, and
$$\lim_{n\to\infty}\frac{\omega(i,I'_{1n}(t),\mathbb H)}{\omega(i,I'_{2n}(t),\mathbb H)}=
\lim_{n\to\infty}\frac{\tan(\pi\omega(i,I'_{1n}(t),\mathbb H))}
{\tan(\pi\omega(i,I'_{2n}(t),\mathbb H))}=$$ $$\lim_{n\to\infty}\frac{\text{meas}\;I'_{1n}(t)
}{\text{meas}\;I'_{2n}(t)}=\lim_{n\to\infty}\frac{\text{meas}\;I_{1n}(t)
}{\text{meas}\;I_{2n}(t)}.$$

This limit exists for every sequence $\{x_n\}$ tending to infinity. So there exists a limit for the
ratio of $m_1(t)$ and $m_2(t)$ as $t\to+0$, and
$$\lim_{t\to+0}\frac{m_1(t)}{m_2(t)}=\lim_{n\to\infty}\frac{m_1({t\over n})}{m_2({t\over n})}=
\lim_{t\to+0}\frac{\text{meas}\;I_1(t)}{\text{meas}\;I_2(t)}=
\frac{\sqrt{c^2+16}+c}{\sqrt{c^2+16}-c}=\frac{1+\beta}{1-\beta},$$ where $\beta$ is given by
(\ref{beta}). This leads to the conclusion desired in Theorem 1 and completes the proof.

\section{Proof of Theorem 2 }

{\it Proof of Theorem 2.} The L\"owner equation (\ref{Low}) admits an explicit integration
\cite{Vas} in the case when $\gamma_0[0,t]$ is a circular arc centered at $i$, $\gamma_0(0)=0$,
with an implicitly given driving function $\lambda(t)$. To be concrete, we will consider
$\gamma_0[0,t]$ such that the argument of $(\gamma_0[0,t]-i)$ increases in $t$. Let a solution
$f_0(z,t)$ to equation (\ref{Low}) map $\mathbb H\setminus\gamma_0[0,t]$ onto $\mathbb H$. Its
inverse $f_0^{-1}(w,t)$ is represented \cite{Vas} by the Christoffel-Schwarz integral
$$\frac{1}{f_0^{-1}(w,t)}=\int_0^{{1\over w}}\frac{(1-\lambda_0w)dw}{(1-\beta_1w)^2(1-\beta_2w)}
=\frac{1}{2\pi}\log\frac{w-\beta_1}{w-\beta_2}+\frac{\beta_2+\beta_1}{\beta_2-\beta_1}\;\frac{1}
{w-\beta_1},$$ where $\beta_1=\beta_1(t)$ and $\beta_2=\beta_2(t)$ are expanded in powers of
$\root3\of t$, $$\beta_1(t)=A_1\root3\of{t^2}+A_2t+A_3\root3\of{t^4}+\dots,,\;\;\;A_1=
-\frac{\root3\of9}{\root3\of{4\pi}},$$ and $$\beta_2(t)=B_1\root3\of
t+B_2\root3\of{t^2}+\dots,\;\;\;B_1=\root3\of{12\pi}.$$ The driving function $\lambda_0(t)$ is
evaluated by $$\lambda_0(t)=2\beta_1(t)+\beta_2(t)=C_1\root3\of
t+C_2\root3\of{t^2}+\dots,\;\;\;C_1=B_1.$$

Suppose that a $C^4$-slit $\gamma_t$ satisfies the conditions of Theorem 2. Then there exists a
driving function $\lambda(t)$ such that a solution $w=f(z,t)$ to equation (\ref{Low}) maps $\mathbb
H\setminus\gamma_t$ onto $\mathbb H$. For the arc-length parameter $s$, represent a transformation
of $\gamma(t(s))$,
\begin{equation}
\tilde\gamma(s):=\frac{2\gamma(t(s))}{2+i\gamma(t(s))}=s+o(s^4),\;\;\;s\to+0. \label{arc1}
\end{equation}

The function $f_0(z,\tau)$ maps $\mathbb H\setminus\gamma_0[0,\tau]$ onto $\mathbb H$. Hence,
$$G_0(w,\tau)=\frac{2f_0^{-1}(w,\tau)}{2+if_0^{-1}(w,\tau)}$$ maps $\mathbb H$ onto the exterior of
the disk of radius 1 centered at $(-i)$ and slit along the segment $[0,\sigma]\subset\mathbb R$.
The length $\sigma(\tau)$ of $[0,\sigma]$ and the half-plane capacity $\tau$ of $\gamma_0[0,\tau]$
are related by
\begin{equation}
\sigma(\tau)=\frac{2f_0^{-1}(\lambda_0(\tau),\tau)}{2+if_0^{-1}(\lambda_0(\tau),\tau)}=
B_1\root3\of{\tau}+O(\root3\of{\tau^2}),\;\;\;\tau\to+0. \label{tau1}
\end{equation}

Let $s$ denote the length of $\tilde\gamma[0,s]$, and let $\sigma$ denote the length of projection
of $\tilde\gamma[0,s]$ onto $[0,\sigma]$ for $\sigma$ large enough. There is a $C^6$-dependence
$s=s(\sigma)$,
$$s(0)=0,\;\;s'(0)=1,\;\;s^{(k)}(0)=0,\;k=2,\dots,6.$$ Therefore,
\begin{equation}
s=\sigma+o(\sigma^6),\;\;\sigma\to+0. \label{sigma1}
\end{equation}

Asymptotic expansion (\ref{arc1}) implies an asymptotic behavior of a distance between
$\tilde\gamma$ and its projection on $[0,\sigma]$,
$$\text{dist}(\tilde\gamma[0,s],[0,\sigma(s)])=\sigma+o(\sigma^6(s)),\;\;s\to+0.$$

According to Lemma 4.10 \cite{Lind}, we have that
$$ t(s)-\tau(\sigma(s))=o(s^3),\;\;s\to+0,$$ where $\sigma(\tau)$ is given by (\ref{tau1}). Hence,
due to (\ref{tau1}) and (\ref{sigma1}),
$$t(s(\sigma))=\tau(\sigma)+o(s^3(\sigma))=\tau(\sigma)+o(\sigma^3)=B_1^{-3}s^3+o(s^3),\;\;
s\to+0.$$ Take into account (\ref{arc1}) and rewrite the last relation in the form
\begin{equation}
\gamma(t)=B_1\root3\of t+\alpha(t)\root3\of t,\;\;\;\lim_{t\to+0}\alpha(t)=0. \label{basic1}
\end{equation}

Choose an arbitrary sequence $\{x_n\}$ of positive numbers $x_n$, $x_n\to\infty$ as $n\to\infty$,
and denote
\begin{equation}
z=g_n(w,t):=\sqrt{x_n}f^{-1}\left(\frac{w}{\sqrt{x_n}},\frac{t}{x_n}\right),\;\;\;n=1,2, \dots\;.
\label{seq1} \end{equation} The function $z=f^{-1}(w,t)$ maps $\mathbb H$ onto $\mathbb
H\setminus\gamma[0,t]$, $f^{-1}(\lambda(t),t)=\gamma(t)$. So the functions
$$G_n(w,t):=\frac{2g_n(w,t)}{2+ig_n(w,t)},\;\;\;n=1,2,\dots,$$ map $\mathbb H$ onto the exterior of
the disk of radius 1 centered at $(-i)$ minus $\tilde\gamma^{(n)}(t)$ where
$$\tilde\gamma^{(n)}(t)=\sqrt{x_n}\,\gamma\left(\frac{t}{x_n}\right)=
B_1\root3\of t+\alpha\left(\frac{t}{x_n}\right)\root3\of t,$$ $$G_n\left(\sqrt{x_n}\;
\lambda\left(\frac{t}{x_n}\right),t\right)=\tilde\gamma^{(n)}(t),\;\;\;0<t\leq T.$$

We see that $$\tilde\gamma^{(n)}(t)-B_1\root3\of t=\alpha\left(\frac{t}{x_n}\right)\root3\of
t\to0,\;\;n\to\infty,$$ and the convergence is uniform with respect to $t\in[0,T]$.

Apply the Rad\'o-Markushevich theorem to $g_n\circ p$ with a conformal mapping $p$ from $\mathbb D$
onto $\mathbb H$ and obtain that the sequence $\{G_n(w,t)\}$ converges to $G_0(w,t)$ which implies
that $\{g_n(w,t)\}$ converges to $f_0^{-1}(w,t)$ as $n\to\infty$ uniformly on compact subsets of
$\mathbb H\cup\mathbb R$.

Denote by $\Gamma_1[0,\tau(t)]$ the {\it left} side of the circular arc $\gamma_0[0,\tau]$ and
denote by $\Gamma_2[0,\tau(t)]$ the {\it right} side of this circular arc. Similarly, denote
$\gamma_{1n}[0,t]$ the {\it left} side of
$$\gamma^{(n)}[0,t]:=\frac{2\tilde\gamma[0,t]}{2-i\tilde\gamma[0,t]},$$ and denote
$\gamma_{2n}[0,t]$ the {\it right} side of $\gamma^{(n)}[0,t]$. The functions $g_n^{-1}(z,t)$ map
$\gamma_{1n}[0,t]$ and $\gamma_{2n}[0,t]$ onto segments $I_{1n}=I_{1n}(t)\subset\mathbb R$ and
$I_{2n}=I_{2n}(t)\subset\mathbb R$, respectively. It is shown \cite{Vas} that slit sides
$\Gamma_1[0,\tau(t)]$ and $\Gamma_2[0,\tau(t)]$ are mapped by $f_0(z,t)$ onto
$$I_1=I_1(t)=[\beta_1(t),\lambda_0(t)]\;\;\text{and}\;\;I_2=I_2(t)=[\lambda_0(t),\beta_2(t)],$$
respectively. The uniform convergence of $g_n$ to $f_0^{-1}$ implies that $I_{1n}(t)$ tend to
$I_1(t)$, and $I_{2n}(t)$ tend to $I_2(t)$ as $n\to\infty$.

Denote by $\gamma'_{1n}[0,t]$ and $\gamma'_{2n}[0,t]$ the {\it left} and the {\it right} sides of
$\gamma[0,{t\over x_n}]$, respectively. The function $f(z,{t\over x_n})$ maps slit sides
$\gamma'_{1n}[0,t]$ and $\gamma'_{2n}[0,t]$ onto segments $I'_{1n}=I'_{1n}(t)\subset\mathbb R$ and
$I'_{2n}=I'_{2n}(t)\subset\mathbb R$, respectively. Compare $I_{kn}(t)$ and $I'_{kn}(t)$ by
(\ref{seq1}) and conclude that $I_{kn}(t)=\sqrt{x_n}I'_{kn}(t)$, and so
$$\text{meas}\;I_{kn}(t)=\sqrt{x_n}\,\text{meas}\;I'_{kn}(t),\;\;k=1,2,\;\;n\geq1,\;\;0<t\leq T.$$

The harmonic measure is invariant under conformal transformations. This gives that $$
\frac{M^2_1({t\over x_n})}{M_2({t\over x_n})}=\frac{\omega^2(f^{-1}(i,{t\over
x_n}),\gamma_1[0,{t\over x_n}],\mathbb H\setminus\gamma({t\over x_n}))}{\omega(f^{-1}(i,{t\over
x_n}),\gamma_2[0,{t\over x_n}],\mathbb H\setminus\gamma({t\over
x_n}))}=\frac{\omega^2(i,I'_{1n}(t),\mathbb H)}{\omega(i,I'_{2n}(t),\mathbb H)}.$$

For $k=1,2$, $n\geq1$, the harmonic measure $\omega(i;I'_{kn}(t),\mathbb H)$ of $I'_{kn}(t)$ at $i$
with respect to $\mathbb H$ equals the angle divided over $\pi$ under which the segment
$I'_{kn}(t)$ is seen from the point $i$. Similarly, the harmonic measure
$\omega(i;I_{kn}(t),\mathbb H)$ of $I_{kn}(t)$ at $i$ with respect to $\mathbb H$ equals the angle
divided over $\pi$ under which the segment $I_{kn}(t)$ is seen from the point $i$. This shows that
the last term in the chain of equalities has a limit as $n\to\infty$, and
$$\lim_{n\to\infty}\frac{\omega^2(i,I'_{1n}(t),\mathbb H)}{\omega(i,I'_{2n}(t),\mathbb H)}=
\lim_{n\to\infty}\frac{\tan^2(\pi\omega(i,I'_{1n}(t),\mathbb H))}
{\tan(\pi\omega(i,I'_{2n}(t),\mathbb H))}=$$ $$\lim_{n\to\infty}\frac{\text{meas}^2\;I'_{1n}(t)
}{\text{meas}\;I'_{2n}(t)}=\lim_{n\to\infty}\frac{\text{meas}^2\;I_{1n}(t)
}{\text{meas}\;I_{2n}(t)}.$$

This limit exists for every sequence $\{x_n\}$ tending to infinity. So there exists a limit for the
ratio of $M^2_1(t)$ and $M_2(t)$ as $t\to+0$, and
$$\lim_{t\to+0}\frac{M^2_1(t)}{M_2(t)}=\lim_{n\to\infty}\frac{M^2_1({t\over n})}{M_2({t\over n})}=
\lim_{t\to+0}\frac{\text{meas}^2\;I_1(t)}{\text{meas}\;I_2(t)}=\lim_{t\to+0}
\frac{(\lambda_0(t)-\beta_1(t))^2}{\beta_2(t)-\lambda_0(t)}=$$ $$\lim_{t\to+0} \frac{(B_1\root3\of
t+(C_2-A_1)\root3\of{t^2}+\dots)^2}{(B_2-C_2)\root3\of{t^2}+\dots}=\frac{B_1^2}{B_2-C_2}=
\frac{B_1^2}{-2A_1}=2\pi.$$ This leads to the conclusion desired in Theorem 2 and completes the
proof.

\end{document}